\documentclass[12pt]{article}
\usepackage{amssymb}

\usepackage{amsmath}
\usepackage{amsmath,amsxtra,amssymb,latexsym, amscd,amsthm}

\advance\voffset-1truecm\relax
\advance\hoffset-1truecm\relax
\def\C{{\Bbb C}}

\newcommand {\Cal}{\mathcal}

\begin{document}

\title{UNIQUENESS PROBLEM\\
FOR MEROMORPHIC MAPPINGS\\
WITH TRUNCATED MULTIPLICITIES\\
AND MOVING TARGETS}
\author{$\quad $Gerd Dethloff and Tran Van Tan}
\date{$\quad$}
\maketitle

\begin{abstract}
\noindent In this paper, using techniques of value distribution theory, we give a
uniqueness theorem for meromorphic mappings of $\mathbb{C}^{m}$ into $%
\mathbb{C}P^{n}$ with $(3n+1)$ moving targets and truncated multiplicities.
\end{abstract}

\section{Introduction}

\indent The uniqueness problem of meromorphic mappings under a condition on
the inverse images of divisors was first studied by R. Nevalinna [6]. He
showed that for two nonconstant meromorphic functions $f$ and $g$ on the
complex plane $\mathbb{C}$, if they have the same inverse images for five
distinct values, then $f \equiv g$, and that $g$ is a special type of linear
fractional tranformation of $f$ if they have the same inverse images,
counted with multiplicities, for four distinct values. These results were
generalized to the case of meromorphic mappings of $\mathbb{C}^m$ into $%
\mathbb{C} P^n$ by H. Fujimoto [1].

In the last years, this problem was continued to be studied by
H.Fujimoto[2], [3], L. Smiley [10], S. Ji [5], M. Ru [9], Z. Tu [11].

Let $f$, $a$ be two meromorphic mappings of $\mathbb{C}^m$ into $\mathbb{C}
P^n$ with reduced representations $f = (f_0 : \dots : f_n)$, $a = (a_0 :
\dots : a_n)$. Set $(f,a):=a_{0}f_{0}+\dots +a_{n}f_{n}$. We say that $a$ is
``small'' with respect to $f$ if $T_{a}(r)=\circ (T_{f}(r))$ as $r\to \infty 
$ (outside a set of finite Lebesgue measure). Assume that $(f,a)%
\not\equiv%
%
0$, denote by $v_{(f,a)}$ the map of $\mathbb{C}^{m}$ into $\mathbb{N}_{0}$
with $v_{(f,a)}(z)=0$ if $(f,a)(z)\neq 0$ and $v_{(f,a)}(z)=k$ if $z$ is a
zero point of $(f,a)$ with multiplicity $k$.

Let $a_{1},\dots ,a_{q}$ \ $(q\geq n+1)$ be meromorphic mappings of $\mathbb{%
C}^{m}$ into $\mathbb{C}P^{n}$ with reduced representations $%
a_{j}=(a_{j0}:\dots :a_{jn})$, \ $j=1,\dots ,q$. We say that $\big\{a_{j}%
\big\}_{j=1}^{q}$ are in general position if for any $1\leq j_{0}<\dots
<j_{n}\leq q$, \ $\text{det}(a_{j_{k}i},0\leq k,i\leq n)%
\not\equiv%
%
0$. 

For each $j\in\{1,...,q\}$, we put $\widetilde{a_j}=(\dfrac{a_{j0}}{a_{jt_j}}:...:\dfrac{a_{jn}}{a_{jt_j}})$ 
and $(f,\widetilde{a_j}) =f_0\dfrac{a_{j0}}{a_{jt_j}}+...+f_n\dfrac{a_{jn}}{a_{jt_j}}$ ,
where $a_{jt_j}$ is the first element of $a_{j0}, ..., a_{jn}$ not identically equal to zero.

Let $%
\mathcal{M}%
%
$ be the field (over $\mathbb{C}$ ) of all meromorphic functions on $\mathbb{%
C}^{m}$ . Denote by $\mathcal{R}\Big(\big\{a_{j}\big\}_{j=1}^{q}\Big)\subset 
\mathcal{M}%
%
$ the subfield generated by the set $\{\dfrac{a_{ji}}{a_{jt_j}},0\leq i\leq n,1\leq j\leq q\}$
over $\mathbb{C}$. Define $\widetilde{\mathcal{R}}\Big(\big\{a_{j}\big\}%
_{j=1}^{q}\Big)\subset 
\mathcal{M}%
%
$ to be the subfield over $\mathbb{C}$ which is generated by all $h\in 
\mathcal{M}%
%
$ with $h^{k}\in \mathcal{R}\Big(\big\{a_{j}\big\}_{j=1}^{q}\Big)$ for some
integer $k$.
 These subfields are independant of the reduced representations
$a_{j}=(a_{j0}:\dots :a_{jn})$, \ $j=1,\dots ,q$, and they are of course also independant of our
choice of the  ${a_{jt_j}}$, because they contain all quotients of the quotients 
$\dfrac{a_{ji}}{a_{jt_j}}, i=0, \dots ,n$.

We say that $f$ is linearly nondegenerate over $\mathcal{R}\Big(\big\{a_{j}%
\big\}_{j=1}^{q}\Big)$ (respectively $\widetilde{\mathcal{R}}\Big(\big\{a_{j}%
\big\}_{j=1}^{q}\Big)$) if $f_{0},\dots ,f_{n}$ are linearly independent
over $\mathcal{R}\Big(\big\{a_{j}\big\}_{j=1}^{q}\Big)$ (respectively $%
\widetilde{\mathcal{R}}\Big(\big\{a_{j}\big\}_{j=1}^{q}\Big)$).

Let $f,g:\mathbb{C}^{m}\longrightarrow \mathbb{C}P^{n}$ be two nonconstant
meromorphic mappings and $\big\{a_{j}\big\}_{j=1}^{q}$ be $q$ ``small''
(with respect to $f$) meromorphic mappings of $\mathbb{C}^{m}$ into $\mathbb{%
\ C}P^{n}$ in general position such that $(f,a_{j}) 
\not\equiv%
%
0$, $(g,a_{j})%
\not\equiv%
%
0$, $j=1,\dots ,q$. Set $\widetilde{E}_{f}^{j}:=\big\{z\in \mathbb{C}
^{m}:v_{(f,a_{j})}(z)>0\big\}$. Assume that:

i) $v_{(f,a_j)} = v_{(g,a_j)}$ for all $j \in \{1, \dots,q\}$

ii) $\text{dim}\big(\widetilde E_f^i \cap \widetilde E_f^j\big)
\leq m-2$ for all $1 \leq i < j \leq q$, and

iii) $f=g$ on $\bigcup\limits_{j=1}^{q}\widetilde{E}_{f}^{j}.$

In [11] Z. Tu showed that:\\

\noindent%
%
\textbf{Theorem A.} If $q=3n+1$ and $f$ is linearly nondegenerate over $%
\mathcal{R}\Big(\big\{
a_{j}\big\}_{j=1}^{q}\Big)$, then there exists a $(n+1)\times (n+1)$-matrix $%
L$ with elements in $\widetilde{\mathcal{R}}\Big(\big\{a_{j}\big\}_{j=1}^{q} %
\Big)$ and $\text{\textrm{det}}(L) 
\not\equiv%
%
0$ such that $L\cdot f=g$.\\
\noindent%
%
\textbf{Theorem B.} If $q=3n+2$ and $f$ is linearly nondegenerate over $%
\widetilde{\mathcal{R}} \Big(\big\{a_{j}\big\}_{j=1}^{q}\Big)$ then $f=g$.%
\\

These theorems (without conditions ii) and iii)) were first showed by H.
Fujimoto ([1]) for hyperplanes ($\big\{a_j\big\}_{j=1}^q$ are constant).

In the above Theorems multiplicities are not truncated (we say that
multiplicities are truncated by a positive integer $M$ if i) is replaced by
the following: $\min \big\{v_{(f,a_{j})},M\big\}=\min \big\{v_{(g,a_{j})},M %
\big\}$). However, the uniqueness problem with truncated multiplicities was
studied in [2], [3], [5], [10] for hyperplanes ($\big\{a_{j}\big\}_{j=1}^{q}$
are constant) and in [9] for moving targets.

For hyperplanes, in [10] L. Smiley proved Theorem B with multiplicities are
truncated by 1, and in [2], [3] H. Fujimoto gave some results related to
Theorem B with multiplicities are truncated by a positive integer $M$.

For moving targets, in [9] M. Ru gave some results related to Theorem B with
multiplicities are truncated by 1, but where the number $q=3n+2$ is replaced
by bigger one.

The main purpose of this paper is to give uniqueness theorems for the case
of $3n+1$ moving targets and multiplicities which are truncated by a
positive integer $M$. Our results are improvements of Theorems A-B where the
number $q=3n+2$ is replaced by smaller one, the multiplicities are truncated
and the condition iii) is replaced by weaker one. In particular, we prove
that for $n \geq 2$ we get $f=g$ already for $q=3n+1$.

The proofs of our results are applications of a generalized Borel Lemma: For
the case where multiplicities are truncated, our object does not satisfy the
assumption ``nowhere vanishing holomorphic functions'' of the (classical)
Borel Lemma. So, first of all, using the techniques of value distribution
theory, we give Lemma 3.1, which is a generalization of the Borel Lemma for
meromorphic functions.

In order to show that under the conditions of our uniqueness theorems the
assumption of Lemma 3.1 is satisfied, we need some results of value
distribution theory of meromorphic mappings of $\mathbb{C}^{m}$ into $%
\mathbb{C}P^{n}$ for moving targets. But the Second Main Theorem as in [8]
(where multiplicities are not truncated) or as in [11] (where multiplicities
are truncated by a positive integer $\ell $) seems to be not sufficient for
our purpose. In order to overcome this difficulty we establish a Second Main
Theorem for meromorphic mappings of $\mathbb{C}^{m}$ into $\mathbb{C}P^{n}$
for $(n+2)$ moving targets with multiplicities truncated by $n$.

Our main results are as follows:\vskip0.25cm

Let $f,g:\mathbb{C}^{m}\longrightarrow \mathbb{C}P^{n}$ be two nonconstant
meromorphic mappings and $\big\{a_{j}\big\}_{j=1}^{3n+1}$ be ``small'' (with
respect to $f$ ) meromorphic mappings of $\mathbb{C}^{m}$ into $\mathbb{C}%
P^{n}$ in general position such that $(f,a_{j})%
\not\equiv%
%
0$, $(g,a_{j})%
\not\equiv%
%
0$, $j=1,\dots ,3n+1$. Put $M=6n(n+1)[N^{2}(N-1)+1]$, where $N=
\begin{pmatrix}
2n+2 \\ 
n+1
\end{pmatrix}
$.

Set $E_{f}^{j}:=\big\{z\in \mathbb{C}^{m}:0\leq v_{(f,a_{j})}(z)\leq M\big\}$%
, ${}^{\ast }E_{f}^{j}:=\big\{z\in \mathbb{C}^{m}:0<v_{(f,a_{j})}(z)\leq M%
\big\}$ and similarly for $E_{g}^{j}$, ${}^{\ast }E_{g}^{j}$, $j=1,\dots
,3n+1$. \\
Assume that:

i) $v_{(f,a_{j})}=v_{(g,a_{j})}$ on $E_{f}^{j}\cap E_{g}^{j}$ for all $j\in
\{1,\dots ,3n+1\}$.

ii) $\text{dim}\big({}^{\ast }E^{i}\cap {}^{\ast }E^{j}\big)\leq m-2$ for
all ${}^{\ast }E^{i}\in \big\{{}^{\ast }E_{f}^{i},{}^{\ast }E_{g}^{i}\big\}$%
, ${}^{\ast }E^{j}\in \big\{{}^{\ast }E_{f}^{j},{}^{\ast }E_{g}^{j}\big\}$
and for all $i\neq j$ with $i\in \{1,\dots ,n+3\},j\in \{1,\dots ,3n+1\}$ .

iii) $f=g$ on $\bigcup\limits_{i=1}^{n+4}\big({}^{\ast }E_{f}^{i}\cap
{}^{\ast }E_{g}^{i}\big)$ for $n\geq 2$ .

This means in particular that in i), ii) and iii) we do not need to pay
attention to points where $v_{(f,a_{j})}$ or $v_{(g,a_{j})}$ is bigger than $M$.%
\\

\noindent%
%
\textbf{Theorem 1.}

\textit{1) If }$n=1$\textit{\ and }$f,g$\textit{\ are linearly nondegenerate
over }$\mathcal{R}\Big(\big\{a_{j}\big\}_{j=1}^{4}\Big)$\textit{\ then there
exists a }$2\times 2$ \textit{-\ matrix}$L$\textit{\ with elements in }$%
\widetilde{\mathcal{R}}\Big(\big\{a_{j}\big\}_{j=1}^{4}\Big)$\textit{\ and }$%
\text{det}(L)%
\not\equiv%
%
0$\textit{\ such that }$L\cdot f=g$\textit{.}

\textit{2) If }$n\geq 2$ and \textit{\ }$f,g$ \textit{\ are linearly
nondegenerate over }$\mathcal{R}\Big(\big\{a_{j}\big\}_{j=1}^{3n+1}\Big)$%
\textit{\ then }$f=g.$\\

We remark that in the case $n=1$, we cannot omit the matrix $L$, as can be
seen easily as follows: Let $f:\mathbb{C} \rightarrow \mathbb{C} $ a
nonconstant nonvanishing holomorphic function, then consider the two functions $f$, $1/f$
and the four values $0, \infty, 1, -1$. Note also that condition i) is
weaker than a truncated multiplicities condition.

We give the following theorem for the case where multiplicities are
truncated.\\

\noindent%
%
\textbf{Theorem 2.} Let $f,g:\mathbb{C}^{m}\longrightarrow \mathbb{C}P^{n}$
be two nonconstant meromorphic mappings and $\big\{a_{j}\big\}_{j=1}^{3n+1}$
be ``small''(with respect to $f$) meromorphic mappings of $\mathbb{C}^{m}$
into $\mathbb{C}P^{n}$ in general position such that $(f,a_{j})%
\not\equiv%
%
0$, $(g,a_{j})%
\not\equiv%
%
0,$ $j=1,\dots ,3n+1.$ Put $M=3n(n+1)N^{2}(N-1)+(3n+4)n$, where $N= 
\begin{pmatrix}
2n+2 \\ 
n+1
\end{pmatrix}
$.\\
Set ${}^{\ast }E_{f}^{j}:=\big\{z\in \mathbb{C}^{m}:0<v_{(f,a_{j})}(z)\leq M%
\big\}$, $j=1,\dots ,3n+1$. \\
Assume that:

\qquad \text{i) min }$\{v_{(f,a_{j})},M\}=$ min $\{v_{(g,a_{j})},M\}$ for
all $j\in \{1,...,3n+1\}.$

\ \ \ \ \ \ ii) dim$\left( ^{\ast }E_{f}^{i}\cap ^{\ast }E_{f}^{j}\right)
\leq m-2$ for all $i\neq j$ with $i\in \{1,...,n+3\},j\in \{1,...,3n+1\}$.

\ \ \ \ \ iii) $f=g$ on $\bigcup\limits_{i=1}^{n+4}$ $^{\ast }E_{f}^{i}$ for 
$n\geq 2.$

Then :

\textit{1) If }$n=1$\textit{\ and }$f$\textit{\ is linearly nondegenerate over 
}$\mathcal{R}\Big(\big\{a_{j}\big\}_{j=1}^{4}\Big)$\textit{\ then there
exists a }$2\times 2$ \textit{-\ matrix}$L$\textit{\ with elements in }$%
\widetilde{\mathcal{R}}\Big(\big\{a_{j}\big\}_{j=1}^{4}\Big)$\textit{\ and }$%
\text{det}(L)%
\not\equiv%
%
0$\textit{\ such that }$L\cdot f=g$\textit{.}

\textit{2) If }$n\geq 2$ \textit{and\ }$f$\textit{\ is linearly
nondegenerate over }$\mathcal{R}\Big(\big\{a_{j}\big\}_{j=1}^{3n+1}\Big)$%
\textit{\ then }$f=g.$\\

We finally remark that in order to obtain uniqueness theorems with \textit{fixed} targets 
only, the authors showed in [13] that the
number $q=3n+1$ of targets can be decreased and that one can use much smaller truncations.\\

\noindent \textbf{Acknowledgements: }The second author would like to thank Professor
Do Duc Thai for valuable discussions, the Universit\'{e} de Bretagne Occidentale
(U.B.O.) for its hospitality and for support, the PICS-CNRS ForMathVietnam
for support.$\qquad $

\section{Preliminaries}

We set $\Vert z\Vert =\big(|z_{1}|^{2}+\dots +|z_{m}|^{2}\big)^{1/2}$ for $%
z=(z_{1},\dots ,z_{m})\in \mathbb{C}^{m}$ and define 
\begin{equation*}
B(r):=\{z\in \mathbb{C}^{m}:|z|<r\},\quad S(r):=\{z\in \mathbb{C}%
^{m}:|z|=r\}\ \text{for all}\ 0<r<\infty .
\end{equation*}
Define \ $d^{c}:=\dfrac{\sqrt{-1}}{4\pi }(\overline{\partial }-\partial )$, $%
\ \upsilon :=(dd^{c}\Vert z\Vert ^{2})^{m-1}$ and 
\begin{equation*}
\sigma :=d^{c}\text{log}\Vert z\Vert ^{2}\wedge (dd^{c}\text{log}\Vert
z\Vert ^{2})^{m-1}.
\end{equation*}
Let $F$ be a nonzero holomorphic function on $\mathbb{C}^{m}$. For each $%
a\in \mathbb{C}^{m}$, expanding $F$ as $F=\sum P_{i}(z-a)$ with homogeneous
polynomials $P_{i}$ of degree $i$ around $a$, we define 
\begin{equation*}
v_{F}(a):=\min \{i:P_{i}%
\not\equiv%
%
0\}.
\end{equation*}
Let $\varphi $ be a nonzero meromorphic function on $\mathbb{C}^{m}$. We
define the map $v_{\varphi }$ as follows: For each $z\in \mathbb{C}^{m}$, we
choose nonzero holomorphic functions $F$ and $G$ on a neighborhood $U$ of $z$
such that $\varphi =\dfrac{F}{G}$ on $U$ and $\text{dim}\big(F^{-1}(0)\cap
G^{-1}(0)\big)\leq m-2$ , and then we put $v_{\varphi }(z):=v_{F}(z)$. Set \ $%
|v_{\varphi }|:=\overline{\big\{z\in \mathbb{C}^{m}:v_{\varphi }(z)\neq 0%
\big\}}$ .

Let $k$, $M$ be positive integers or $+\infty $. Set 
\begin{equation*}
{}^{\leq M}v_{\varphi }^{[k]}(z)=0\text{ if}\quad v_{\varphi }(z)>M\text{
and }^{\leq M}v_{\varphi }^{[k]}(z)=\min \{v_{\varphi }(z),k\}\text{ }\text{%
if}\quad v_{\varphi }(z)\leq M
\end{equation*}
\begin{align*}
{}^{>M}v_{\varphi }^{[k]}(z)=0\text{ }& \text{if}\quad v_{\varphi }(z)\leq M%
\text{ and }^{>M}v_{\varphi }^{[k]}(z)=\min \{v_{\varphi }(z),k\}\text{ if }%
v_{\varphi }(z)>M. \\
& {}
\end{align*}
We define 
\begin{equation*}
{}^{\leq M}N_{\varphi }^{[k]}(r):=\int\limits_{1}^{r}\frac{{}^{\leq M}n(t)}{%
t^{2m-1}}dt
\end{equation*}
and 
\begin{equation*}
{}^{>M}N_{\varphi }^{[k]}(r):=\int\limits_{1}^{r}\frac{{}^{>M}n(t)}{t^{2m-1}}%
dt\qquad (1\leq r<+\infty )
\end{equation*}
where 
\begin{align*}
{}^{\leq M}n(t):={\int\limits_{|v_{\varphi }|\cap B(r)}}{}^{\leq
M}v_{\varphi }^{[k]}.\upsilon \quad & \text{for}\quad m\geq 2{}\text{ , }%
^{\leq M}n(t):={\sum\limits_{|z|\leq t}}{}^{\leq M}v_{\varphi
}^{[k]}(z)\quad \text{for}\quad m=1 \\
{}^{>M}n(t):={\int\limits_{|v_{\varphi }|\cap B(r)}}{}^{>M}v_{\varphi
}^{[k]}.\upsilon \quad & \text{for}\quad m\geq 2\text{ , }{}^{>M}n(t):={%
\sum\limits_{|z|\leq t}}{}^{>M}v_{\varphi }^{[k]}(z)\quad \text{for}\quad
m=1.
\end{align*}

Set \ $N_{\varphi }(r):={}^{\leq \infty }N_{\varphi }^{[\infty ]}(r)$, \ $%
N_{\varphi }^{[k]}(r):={}^{\leq \infty }N_{\varphi }^{[k]}(r)$.

We have the following Jensen's formula (see [3], P.177): 
\begin{equation*}
N_{\varphi }(r)-N_{\frac{1}{\varphi }}(r)=\int\limits_{S(r)}\text{log}%
|\varphi |\sigma -\int\limits_{S(1)}\text{log}|\varphi |\sigma .
\end{equation*}

Let $f : \mathbb{C}^m \longrightarrow \mathbb{C} P^n$ be a meromorphic
mapping. For arbitrary fixed homogeneous coordinates $(w_0 : \dots : w_n)$
of $\mathbb{C} P^n$, we take a reduced representation $f = (f_0 : \dots :
f_n)$ , which means that each $f_i$ is a holomorphic function on $\mathbb{C}^n$
and $f(z) = (f_0(z) : \dots : f_n(z))$ outside the analytic set $\{f_0 =
\dots = f_n = 0\}$ of codimension $\geq 2$. Set $\Vert f\Vert = \big(|f_0|^2
+ \dots + |f_n|^2\big)^{1/2}$.

The characteristic function of $f$ is defined by 
\begin{equation*}
T_{f}(r)=\int\limits_{S(r)}\text{log}\Vert f\Vert \sigma -\int\limits_{S(1)}%
\text{log}\Vert f\Vert \sigma ,\quad 1\leq r<+\infty .
\end{equation*}
For a meromorphic function $\varphi $ on $\mathbb{C}^{m}$, the proximity function is defined by
$$m(r,\varphi): = \int\limits_{S(r)}\log^+|\varphi|\sigma$$
 and we have, by the classical First Main Theorem that (see [4], p.135) 
 $$m(r,\varphi) \leq T_{\varphi }(r) +O(1).$$
 Here,
 the characteristic function $T_{\varphi }(r)$ of $\varphi $ is defined \ by considering $%
\varphi $ as a meromorphic mapping of $\mathbb{C}^{m}$ into $\mathbb{C}%
P^{1}. $

We state the First and Second Main Theorem of Value Distribution Theory (see
e.g.[11], [2]):

\vskip0.25cm \noindent \textbf{First Main Theorem.} (Moving target version)
\textit{\ Let $a$ be a
meromorphic mapping of  $\mathbb{C}^{m}$ into $\mathbb{C}P^{n}$ such that $%
(f,a)%
\not\equiv%
%
0$. Then for reduced representations   $f=(f_{0}:\dots :f_{n})$ and $%
a=(a_{0}:\dots :a_{n})$, we have:
\begin{equation*}
N_{(f,a)}(r)\leq T_{f}(r)+T_{a}(r)\quad \text{for all}\quad r\geq 1.
\end{equation*}
}

For a hyperplane $H:a_{0}w_{0}+\dots +a_{n}w_{n}=0$ in $\mathbb{C}P^{n}$
with $\text{im}\,f\nsubseteq H$, we denote $(f,H)=a_{0}f_{0}+\dots
+a_{n}f_{n}$ , where $(f_{0}:\dots :f_{n})$ again is a reduced
representation of $f$.

\vskip0.25cm \noindent \textbf{Second Main Theorem.} (Classical version) \textit{Let $f$ be a
linearly nondegenerate meromorphic mapping of $\mathbb{C}^{m}$ into $\mathbb{%
C}P^{n}$ and $H_{1},\dots ,H_{q}$ $(q\geq n+1)$ hyperplanes of $\mathbb{C}%
P^{n}$ in general position, then 
\begin{equation*}
(q-n-1)T_{f}(r)\leq
\sum\limits_{j=1}^{q}N_{(f,H_{j})}^{[n]}(r)+o(T_{f}(r))\quad
\end{equation*}
}for all$\ r$ \textit{except for a set of finite Lebesgue measure.}

\section{Proof of our results}

First of all, we give a generalization of the Borel Lemma for meromorphic
functions.\\

\noindent%
%
\textbf{Lemma 3.1.} Let $h_{0},\dots ,h_{t}$ $(t\geq 2)$ be nonzero
meromorphic functions on $\mathbb{C}^{m}$ and $A$ be a subset of $\
(1,+\infty )$ with infinite Lebesgue measure. Assume that

\text{a)} $h_{0}+\dots +h_{t}\equiv 0,$

\vskip0.15cm \text{b)} $\sum\limits_{v=0}^{t}N_{h_{v}}^{[1]}(r)+\sum%
\limits_{v=0}^{t}N_{\frac{1}{h_{v}}}^{[1]}(r)\leq \dfrac{1}{t(t+1)}%
T_{\varphi _{ijk}}(r)$, \ $r\in A$ for all
$\{i,j,k\}\subset \{0,1,\dots ,t\}$ such that $\dfrac{h_{i}}{h_{j}}$%
\thinspace , $\dfrac{h_{j}}{h_{k}}$\thinspace , $\dfrac{h_{k}}{h_{i}}$ are
all nonconstant, where $\varphi _{ijk}:=[h_{i}:h_{j}:h_{k}]$ is a
meromorphic mapping of $\mathbb{C}^{m}$ into $\mathbb{C}P^{2}$.

Then there exists a decomposition of indices $\{0,\dots ,t\}=I_{1}\cup \dots
\cup I_{s}$ such that:

\text{i)} $\#I_{v}\geq 2$ for all $v\in \{1,\dots ,s\},$

\text{ii)} $i,j\in I_{v}$ if and only if $\dfrac{h_{i}}{h_{j}}$ is constant,

\vskip0.15cm \text{iii)} $\sum\limits_{j\in I_{v}}h_{j}=0$, \ $v\in
\{1,\dots ,s\}$.\\

\noindent%
%
\textbf{Proof.}
We prove this Lemma by induction on $t.$

+) If $t=2$, we have 
\begin{equation}
h_{0}+h_{1}+h_{2}\equiv 0 .  \tag{1}
\end{equation}

\noindent \textbf{Case 1.} If one of the meromorphic functions $\dfrac{h_{0}%
}{h_{1}}$\thinspace , $\dfrac{h_{1}}{h_{2}}$\thinspace , $\dfrac{h_{2}}{h_{0}%
}$ is constant, then by (1), we have that $h_{0}:h_{1}:h_{2}$ are
constant. We get i), ii) and iii).

\vskip0.15cm \noindent \textbf{Case 2.} If $\dfrac{h_{0}}{h_{1}}$\thinspace
, $\dfrac{h_{1}}{h_{2}}$\thinspace , $\dfrac{h_{2}}{h_{0}}$ are nonconstant,
by Theorem 5.2.29 in [7], we have 
\begin{equation*}
T_{\varphi _{012}}(r)=T_{[h_{0}:h_{1}:h_{2}]}(r)\leq
T_{[h_{0}:h_{1}]}(r)+T_{[h_{0}:h_{2}]}(r)+0(1).
\end{equation*}
Without loss of generality, we may assume that $T_{[h_{0}:h_{1}]}(r)\geq
T_{[h_{0}:h_{2}]}(r)$ for all $r\in A_{1}$, where $A_{1}$ is a subset of $A$
with infinite Lebesgue measure. Then 
\begin{equation*}
\sum_{i=0}^{2}N_{h_{i}}^{[1]}(r)+\sum_{i=0}^{2}N_{\frac{1}{h_{i}}%
}^{[1]}(r)\leq \frac{1}{6}T_{\varphi _{012}}(r)\leq \frac{1}{3}%
T_{[h_{0}:h_{1}]}(r),\quad r\in A_{1}.
\end{equation*}
Let $[h_{0}^{\prime }:h_{1}^{\prime }]$ be a reduced representation of $%
[h_{0}:h_{1}]:\mathbb{C}^{m}\longrightarrow \mathbb{C}P^{1}$, $h_{0}^{\prime
}$ and $h_{1}^{\prime }$ are holomorphic functions. Set $h_{2}^{\prime }=%
\dfrac{h_{0}^{\prime }h_{2}}{h_{0}}$, then 
\begin{equation*}
h_{0}^{\prime }+h_{1}^{\prime }+h_{2}^{\prime }\equiv 0.
\end{equation*}
For each $j\in \{0,1,2\}$, we have that a zero of $h_{j}^{\prime }$ is a pole
or a zero of some $h_{i}$ $(i\in \{0,1,2\})$. On the other hand 
\begin{equation*}
\text{dim}\{z:h_{0}^{\prime }(z)=h_{1}^{\prime }(z)=0\}\leq m-2.
\end{equation*}
Hence, we get
\begin{equation*}
\sum_{i=0}^{2}N_{h_{i}^{\prime }}^{[1]}(r)\leq 2\cdot \Big(%
\sum_{i=0}^{2}N_{h_{i}}^{[1]}(r)+\sum_{i=0}^{2}N_{\frac{1}{h_{i}}}^{[1]}(r)%
\Big)\leq \frac{2}{3}T_{[h_{0}:h_{1}]}(r),\quad r\in A_{1}.
\end{equation*}
By the Second Main Theorem, we have: 
\begin{align*}
T_{[h_{0}:h_{1}]}(r)& \leq N_{h_{0}^{\prime }}^{[1]}(r)+N_{h_{1}^{\prime
}}^{[1]}(r)+N_{h_{0}^{\prime }+h_{1}^{\prime }}^{[1]}(r)+\circ \big(%
T_{[h_{0}:h_{1}]}(r)\big) \\
& =\sum_{i=0}^{2}N_{h_{i}^{\prime }}^{[1]}(r)+\circ \big(T_{[h_{0}:h_{1}]}(r)%
\big) \\
& \leq \frac{2}{3}T_{[h_{0}:h_{1}]}(r)+\circ \big(T_{[h_{0}:h_{1}]}(r)\big)%
,\quad r\in A_{1}.
\end{align*}
This is a contradiction when $r\rightarrow \infty $, $r\in A_{1}$.

This completes the proof of the case $t=2$.\\

+) Assume that our assertion holds up to $t$ $(t\geq 2)$. Consider 
\begin{equation}
h_{0}+\dots +h_{t+1}\equiv 0.  \tag{2}
\end{equation}
We introduce an equivalence relation in $\{0,\dots ,t+1\}$ as follows: $%
i\sim j$ if and only if $\dfrac{h_{i}}{h_{j}}$ is constant. Let 
\begin{equation*}
\{I_{1},\dots ,I_{s}\}=\{0,\dots ,t+1\}/\sim .
\end{equation*}
By definition we have ii).

For the proof of i), we assume that there exists $I_{v}$ containing only one
index, say $I_{s}=\{t+1\}$. Then $\dfrac{h_{i}}{h_{t+1}}$ $(i=0,\dots ,t)$
are all nonconstant.

If $s=2$ then $I_{1}=\{0,\dots ,t\}$, $I_{2}=\{t+1\}$.

By (2) we have 
\begin{equation*}
c\cdot h_{0}+h_{t+1}\equiv 0,\quad c\in \mathbb{C}^{\ast }.
\end{equation*}

Thus $\dfrac{h_{0}}{h_{t+1}}$ is constant, this is a contradiction.

If $s=3$, without loss of generality we may assume that $0\in I_{1}$, $1\in
I_{2}$. \\
By (2) we have
\begin{equation*}
c\cdot h_{0}+d\cdot h_{1}+h_{t+1}\equiv 0,\quad c,d\in \mathbb{C}\text{ .}
\end{equation*}

* If $c\cdot d=0$, then $t+1\in I_{1}$ or $t+1\in I_{2}$, this is a
contradiction.

* If $c\neq 0$, $d\neq 0$, we have: 
\begin{equation*}
T_{[c\cdot h_{0}:d\cdot h_{1}:h_{t+1}]}(r)=T_{[h_{0}:h_{1}:h_{t+1}]}(r)+0(1)%
\text{ .}
\end{equation*}
So by the basic step of induction, we have that $h_{0}:h_{1}:h_{t+1}$ are
constant. This is a contradiction.

If $s>3$, let $\Psi :=[h_{0}:\dots :h_{t}]:\mathbb{C}^{m}\longrightarrow 
\mathbb{C}P^{t}$.

Let $[h_{0}^{\prime }:\dots :h_{t}^{\prime }]$ be a reduced representation
of $\Psi $.

\vskip0.15cm Set $h_{t+1}^{\prime }=\dfrac{h_{t}^{\prime }\cdot h_{t+1}}{%
h_{t}}$\thinspace , then $h_{0}^{\prime }+\dots +h_{t+1}^{\prime }\equiv 0$.

\vskip0.15cm For each $j\in \{0,\dots ,t+1\}$, we have that a zero of $%
h_{j}^{\prime }$ is a pole or a zero of some $h_{i}$ $(i\in \{0,\dots
,t+1\}) $ .

Hence, we get
\begin{align*}
N_{h_{j}^{\prime }}^{[1]}(r)& \leq
\sum_{i=0}^{t+1}N_{h_{i}}^{[1]}(r)+\sum_{i=0}^{t+1}N_{\frac{1}{h_{i}}%
}^{[1]}(r) \\
& \leq \frac{1}{(t+1)(t+2)}T_{\varphi _{kpq}}(r) \\
& \leq \frac{1}{(t+1)(t+2)}T_{\Psi }(r),\quad r\in A
\end{align*}
where $k\in I_{1}$, $p\in I_{2}$, $q\in I_{3}$.

If $\Psi $ is linearly nondegenerate, by the Second Main Theorem we have: 
\begin{align*}
T_{\Psi }(r)& \leq \sum_{i=0}^{t}N_{h_{i}^{\prime
}}^{[t]}(r)+N_{h_{0}^{\prime }+\dots +h_{t}^{\prime }}^{[t]}(r)+\circ \big(%
T_{\Psi }(r)\big) \\
& =\sum_{t=0}^{t+1}N_{h_{i}^{\prime }}^{[t]}(r)+\circ \big(T_{\Psi }(r)\big)%
\leq t\cdot \sum_{i=0}^{t+1}N_{h_{i}^{\prime }}^{[1]}(r)+\circ \big(T_{\Psi
}(r)\big) \\
& \leq \frac{t(t+2)}{(t+1)(t+2)}T_{\Psi }(r)+\circ \big(T_{\Psi }(r)\big) \\
& =\frac{t}{t+1}T_{\Psi }(r)+\circ \big(T_{\Psi }(r)\big),\quad r\in A\text{
.}
\end{align*}

This is a contradiction when $r\rightarrow \infty $, $r\in A$.

Thus, $\Psi $ is linearly degenerate, so there exist constants

$(C_{0},\dots ,C_{t})\neq (0,\dots ,0)$ such that 
\begin{equation}
C_{0}h_{0}+\dots +C_{t}h_{t}=0\text{ .}  \tag{3}
\end{equation}
We may assume that $C_{0}=1$. By (2) and (3) we have 
\begin{equation*}
(C_{1}-1)h_{1}+\dots +(C_{t}-1)h_{t}-h_{t+1}\equiv 0\text{ .}
\end{equation*}
It can be written in the form: 
\begin{equation}
a_{1}h_{i_{1}}+\dots +a_{k}h_{i_{k}}+a_{t+1}h_{t+1}\equiv 0  \tag{4}
\end{equation}
such that $a_{i}\in \mathbb{C}^{\ast }$, $a_{t+1}=-1$, $\dfrac{h_{p}}{h_{q}}$
is nonconstant for all $p\neq q\in \{i_{1},\dots ,i_{k},t+1\}$ and $k\leq
t-1 $.

+) If $k=1$, by (4) we have that $h_{i_{1}}:h_{t+1}$ is constant. This is a
contradiction.

+) If $k\geq 2$, for $\{p,q,v\}\subset \{i_{1},\dots ,i_{k},t+1\}$ we have 
\begin{equation*}
T_{[a_{p}h_{p}:a_{q}h_{q}:a_{v}h_{v}]}(r)=T_{[h_{p}:h_{q}:h_{v}]}(r)+0(1).
\end{equation*}

By the induction hypothesis (since $k+1\leq t$) there exists $p\in
\{i_{1},\dots ,i_{k}\}$ such that $a_{p}h_{p}:a_{t+1}h_{t+1}$ is constant.
Thus $h_{p}:h_{t+1}$ is constant, this is a contradiction.

So $\#I_{v}\geq 2$ for all $v\in \{1,\dots ,s\}$, we get i).

Finally we show iii). We choose an index $v\in I_{v}$ and set 
\begin{equation*}
\sum_{i\in I_{v}}h_{i}=c_{v}\cdot h_{v},\quad c_{v}\in \mathbb{C}\text{ .}
\end{equation*}
Then (2) can be written as 
\begin{equation*}
\sum_{v=1}^{s}c_{v}\cdot h_{v}\equiv 0
\end{equation*}
By i) and the induction hypothesis, we infer like above that $c_{v}\equiv 0$%
. This shows iii). We have completed the proof of Lemma 3.1. $\square $%
\\

We give the Second Main Theorem of meromorphic mappings of $\mathbb{C}^{m}$
into \textit{$\mathbb{C}P^{n}$ }with $(n+2)$ moving targets.\\

\noindent%
%
\textbf{Lemma 3.2.} \textit{Let }$f,g:\mathbb{C}^{m}\longrightarrow \mathbb{C}%
P^{n}$\textit{\ be nonconstant meromorphic mappings and }$\big\{a_{j}\big\}%
_{j=1}^{n+2}$\textit{\ be ``small'' (with respect to} $g$ )\textit{\
meromorphic mappings of }$\mathbb{C}^{m}$\textit{\ into }$\mathbb{C}P^{n}$%
\textit{\ in general position.} 

a)
Denote the meromorphic mapping,
\begin{align*}
F = \big( c_1 \cdot (f,\tilde a_1) : \cdots :
c_{n+1} \cdot (f,\tilde a_{n+1})\big) : \C^m \longrightarrow \C P^n
\end{align*}
where $\big\{c_i\big\}_{i=1}^{n+1}$ are ``small" (with respect to $g$)
nonzero meromorphic functions on $\C^m$.
Then we have
\begin{align*}
T_F(r) = T_f(r) + o(T_g(r)).
\end{align*}
Moreover, if 
\begin{align*}
f &= (f_1: \cdots : f_{n+1}),\\
a_i &= (a_{i1} : \cdots : a_{i(n+1)}),\\
F &= \Big( \frac{c_1 \cdot (f,\tilde a_1)}{h} : \cdots :
\frac{c_{n+1} \cdot (f,\tilde a_{n+1})}{h}\Big)
\end{align*}
are reduced representations, where $h$ is a meromorphic function
on $\C^m$, then
\begin{align*}
N_{h}(r) \leq o(T_g(r))
\end{align*}
and
\begin{align*}
N_{\frac{1}{h}}(r) \leq o(T_g(r)).
\end{align*}

b)
\textit{ Assume that} $f$ \textit{is linearly nondegenerate over}
$\mathcal{R}\Big(\big\{a_{j}\big\}_{j=1}^{n+2}\Big).$
\textit{Then we have} 
\begin{equation*}
T_{f}(r)\leq \sum\limits_{j=1}^{n+2}N_{(f,a_{j})}^{[n]}(r)+o(T_{f}(r))+o(T_g(r))\quad 
\end{equation*}
\textit{for all}$\ r$ \textit{except for a set of finite Lebesgue measure.}
\\

\noindent \textbf{Proof.} 

a) Set
\begin{align*}
F_i = \frac{c_i \cdot (f,\tilde a_i)}{h}\,,\quad
i \in \{1,\dots,n+1\}.
\end{align*}
So we have
\begin{equation}
\begin{cases}
c_1 a_{10}f_0 + \cdots + c_1 a_{1n}f_n = h \cdot F_1 \cdot a_{1t_1}\\
\dots \quad\qquad \dots \qquad \quad \dots \qquad \quad \dots \\
c_{n+1}a_{(n+1)0}f_0 + \cdots + c_{n+1}a_{(n+1)n}f_n =
h \cdot F_{n+1} \cdot a_{(n+1)t_{n+1}}
\end{cases} \tag{5}
\end{equation}
Since $(F_1 : \cdots : F_{n+1})$ is a reduced representation of $F$,
$\text{codim}\,\{ F_1 = \cdots = F_{n+1} = 0\} \geq 2$.
Hence, by (5) we see:
\begin{align*}
N_{\frac{1}{h}} (r) \leq \sum_{i=1}^{n+1} N_{a_{i t_i}}(r) +
\sum_{i=1}^{n+1} N_{\frac{1}{c_i}}(r) = o(T_g(r)).
\end{align*}
Set
\begin{align*}
P := \begin{pmatrix}
c_1 a_{10} & \dots &c_1 a_{1n}\cr
\vdots & \ddots & \vdots \cr
c_{n+1}a_{(n+1)0}& \dots & c_{n+1}a_{(n+1)n}\end{pmatrix} ,
\end{align*}
and matrices $P_i$ $(i \in \{1,\dots,n+1\})$ which are
defined from $P$ after changing the $i^{th}$ column by
$\begin{pmatrix} F_1 a_{1t_1}\cr
\vdots \cr
F_{(n+1)}a_{(n+1)t_{n+1}}\end{pmatrix}$.

Put $u = \text{det}(P)$ and $u_i = \text{det}(P_i)$, $i \in \{1, \dots, n+1\}$.
It is easy to see that:
\begin{align*}
\frac{u}{a_{1t_1} \cdots a_{(n+1)t_{n+1}}} \in 
\Cal R\big(\big\{a_j\big\}_{j=1}^{n+2}\big)
\end{align*}
and
\begin{align*}
N_{\frac{1}{u_i}}(r) \leq 
0\Big(\sum_{j=1}^{n+1} N_{\frac{1}{c_j}}(r)\Big) = 
o(T_g(r)),\quad i = 1,\dots, n+1.
\end{align*}
By (5) we have
\begin{align}
\begin{cases} f_0 = \dfrac{h \cdot u_1}{u} \cr
\quad \quad \dots \cr
f_n = \dfrac{h \cdot u_{n+1}}{u} \end{cases} \tag{6}
\end{align}
On the other hand $(f_0 : \cdots : f_n)$ is a reduced
representation of $f$. Hence,

\begin{align*}
N_h(r) &\leq N_u(r) + \sum_{i=1}^{n+1}
N_{\frac{1}{u_i}}(r)\\
&\leq N_{\frac{u}{a_{1t_1} \cdots a_{(n+1)t_{n+1}}}}(r)
+ N_{a_{1t_1}\cdots a_{(n+1)t_{n+1}}}(r) 
+ \sum_{i=1}^{n+1} N_{\frac{1}{u_i}}(r)\\
&= o(T_g(r)).
\end{align*}
We have

\begin{align}
T_F(r) &= \int\limits_{S(r)} \text{log} \Big(\sum_{i=1}^{n+1}
|F_i|^2\Big)^{1/2} \sigma + 0(1)\notag\\
&= \int\limits_{S(r)} \text{log}\Big(\sum_{i=1}^{n+1}
\Big|\frac{c_i \cdot (f, \tilde a_i)}{h}\Big|^2 \Big)^{1/2}
\sigma + 0(1)\notag\\
&= \int\limits_{S(r)} \text{log}\Big( \sum_{i=1}^{n+1}
|c_i(f,\tilde a_i)|^2\Big)^{1/2} \sigma
- \int\limits_{S(r)} \text{log}|h| \sigma + 0(1)\notag\\
&\leq \int\limits_{S(r)} \text{log}\Vert f \Vert \sigma +
\int\limits_{S(r)} \text{log}\Bigg(\sum_{i=1}^{n+1}|c_i|^2
\Big(\Big|\frac{a_{i0}}{a_{it_i}}\Big|^2 + \cdots +
\Big|\frac{a_{in}}{a_{it_i}}\Big|^2\Big)\Bigg)^{1/2} \sigma \notag\\
&\quad - N_h(r) + N_{\frac{1}{h}}(r) + 0(1)\notag\\
&\leq T_f(r) + \int\limits_{S(r)}
\text{log}^+\Bigg(\sum_{i=1}^{n+1}
\Big(\Big|c_i \cdot \frac{a_{i0}}{a_{it_i}}\Big|^2 + \cdots +
\Big|c_i \cdot \frac{a_{in}}{a_{it_i}}\Big|^2\Big)\Bigg)^{1/2}\sigma
+ o(T_g(r))\notag\\
&\leq T_f(r) + \sum_{i=1}^{n+1} \sum_{j=0}^n 
m\Big(r,c_i \frac{a_{ij}}{a_{it_i}}\Big)
+ o(T_g(r))\notag\\
&= T_f(r) + o(T_g(r)) \tag{7}
\end{align}
(note that $c_i \cdot \dfrac{a_{ij}}{a_{it_i}} \in
\Cal R\big(\big\{a_j\big\}_{j=1}^{n+2}\big)$).

\noindent
(6) can be written as
\begin{align*}
\begin{cases}
f_0 &= h \cdot \sum\limits_{j=1}^{n+1} b_{j0}F_j\\
\dots &\quad \dots \quad \dots \\
f_n &= h \cdot \sum\limits_{j=1}^{n+1} b_{jn} F_j \end{cases}
\end{align*}
where \ $b_{ji} \in \Cal R \big(\big\{a_j\big\}_{j=1}^{n+2} \big)$.

So we get
\begin{align}
T_f(r) &= \int\limits_{S(r)} \text{log}\Vert f\Vert \sigma + 
0(1) \notag\\
&= \int\limits_{S(r)} \text{log} \Big(\sum_{i=0}^n
\Big| \sum_{j=1}^{n+1} b_{ji}F_j \Big|^2\Big)^{1/2} \sigma
+ \int\limits_{S(r)}\text{log}|h|\sigma + 0(1)\notag\\
&\leq \int\limits_{S(r)} \text{log}\Vert F \Vert \sigma +
\int\limits_{S(r)} \text{log} \Big( \sum_{i, j} |b_{ji}|^2
\Big)^{1/2} \sigma + N_h(r) - N_{\frac{1}{h}}(r)
+ 0(1) \notag\\
&\leq T_F(r) + \int\limits_{S(r)} \text{log}^+
\Big(\sum_{i,j}|b_{ji}|^2\Big)^{1/2}\sigma + o(T_g(r))\notag\\
&\leq T_F(r) + \sum_{i,j}m(r,b_{ij}) + o(T_g(r)) \notag\\
&= T_F(r) + o(T_g(r)). \tag{8}
\end{align}
By (7) and (8), we have
\begin{align*}
T_F(r) = T_f(r) + o(T_g(r)) .
\end{align*}
This finishes the proof of part a).\\

b) 
We use a) for a special set of $c_i$ :
Set
\begin{align*}
N_{n+2} := \begin{pmatrix}
a_{10}& \dots &a_{(n+1)0}\cr
a_{11}& \dots &a_{(n+1)1}\cr
\vdots& \ddots& \vdots \cr
a_{1n}& \dots & a_{(n+1)n} \end{pmatrix}
\end{align*}
and matrices $N_i$, $i \in \{1,\dots,n+1\}$, which are defined from
$N_{n+2}$ after changing the $i^{th}$ column by
$\begin{pmatrix} a_{(n+2)0}\cr
\vdots \cr
a_{(n+2)n}\end{pmatrix}$.

Set
\begin{align*}
c_i = \frac{a_{i t_i}}{a_{1 t_1} \cdots a_{(n+2)t_{n+2}}}\cdot
\text{det}(N_i),\quad i \in \{1,\dots,n+2\},
\end{align*}
then
\begin{align*}
c_i \in \Cal R\big(\big\{a_j\big\}_{j=1}^{n+2}\big),\quad 
i \in \{1, \dots, n+2\}.
\end{align*}
It is easy to see that:
\begin{align}
\sum_{i=1}^{n+1} c_i \cdot (f, \tilde a_i) =
c_{n+2} \cdot (f, \tilde a_{n+2}) \tag{9}
\end{align}

$F$ is a linearly nondegenerate meromorphic mapping,
since $f$ is linearly nondegenerate over 
$\Cal R\big(\big\{a_j\big\}_{j=1}^{n+2}\big)$
and since the $a_j$ $(j = 1,\dots,n+2)$ are in general position.

Thus, by the First and the Second Main Theorem, we have

\begin{align*}
T_f(r)+o(T_g(r))
&= T_F(r) \leq \sum_{j=1}^{n+1}N_{F_j}^{[n]}(r) +
N_{F_1+\cdots + F_{n+1}}^{[n]}(r) + o(T_F(r))\\
&\overset {(9)}{=} \sum_{j=1}^{n+2}
N_{\frac{c_j\cdot (f,\tilde a_j)}{h}}^{[n]}(r) +
o(T_f(r))+o(T_g(r))\\
&\leq \sum_{j=1}^{n+2}N_{(f,\tilde a_j)}^{[n]}(r) +
\sum_{j=1}^{n+2}N_{c_j}(r) + 
(n+2)N_{\frac{1}{h}}(r) + o(T_f(r))+o(T_g(r))\\
&= \sum_{j=1}^{n+2} N_{(f,\tilde a_j)}^{[n]}(r) + o(T_f(r)) +o(T_g(r))\\
&= \sum_{j=1}^{n+2} N_{(f, a_j)}^{[n]}(r)  + o(T_f(r)) +o(T_g(r)).
\end{align*}
This completes proof of Lemma 3.2. $\square $\\

\noindent%
%
\textbf{Proof of Theorem 1.}
Without loss of generality, we may assume that there exists a subset $A$ of $%
(1,+\infty )$ with infinite Lebesgue measure such that 
\begin{equation}
T_{f}(r)\geq T_{g}(r),\quad r\in A,  \tag{10}
\end{equation}
(note that if $T_{g}(r)\geq T_{f}(r)$ for all $r$ except for a set of finite
Lebesgue measure then $\big\{a_{j}\big\}_{j=1}^{3n+1}$ are ``small'' with
respect to $g$). Define functions 
\begin{equation*}
h_{j}=\frac{(f,a_{j})}{(g,a_{j})}\,,\quad j\in \{1,\dots ,3n+1\}.
\end{equation*}

We choose an arbitrary subset $Q=\{j_{1},\dots ,j_{2n+2}\}$ of the index set 
$Q_{0}:=\{1,\dots ,3n+1\}$.

We now prove that:

For each $I\subset Q$, $\#I=n+1$, there exists some $J\subset Q$ with $I\neq
J$, $\#J=n+1$ such that 
\begin{equation}
\dfrac{h_{I}}{h_{J}}\in \mathcal{R}\Big(\big\{a_{j}\big\}_{j=1}^{3n+1}\Big)%
,\quad \text{where}\quad h_{I}=\prod\limits_{i\in I}h_{i}.  \tag{11}
\end{equation}
We have 
\begin{equation*}
\left\{ 
\begin{array}{c}
a_{j0}f_{0}+\dots +a_{jn}f_{n}=h_{j}(a_{j0}g_{0}+\dots +a_{jn}g_{n}) \\ 
j\in Q
\end{array}
\right.
\end{equation*}
\begin{equation*}
\Rightarrow \left\{ 
\begin{array}{c}
a_{j_{s}0}f_{0}+\dots +a_{j_{s}n}f_{n}-h_{j_{s}}a_{j_{s}0}g_{0}-\dots
-h_{j_{s}}a_{j_{s}n}g_{n}=0 \\ 
1\leq s\leq 2n+2
\end{array}
\right.
\end{equation*}
Therefore, we get 
\begin{equation*}
\text{det}(a_{j_{s}0},\dots ,a_{j_{s}n},h_{j_{s}}a_{j_{s}0},\dots
,h_{j_{s}}a_{j_{s}n},\ 1\leq s\leq 2n+2)\equiv 0\text{ .}
\end{equation*}
For each $I=\{j_{s_{0}},\dots ,j_{s_{n}}\}\subset Q$, $1\leq s_{0}<\dots
<s_{n}\leq 2n+2$, we define 
\begin{equation*}
A_{I}=\frac{(-1)^{\frac{n(n+1)}{2}+s_{0}+\dots +s_{n}}\cdot \text{det}%
(a_{j_{s_{k}}i},0\leq k,i\leq n)\cdot \text{det}(a_{j_{s_{k}^{\prime
}}i},0\leq k,i\leq n)}{a_ {j_1t_{j_1}} \dots a_ {j_{2n+2}t_{j_{2n+2}}}}
\end{equation*}
where $\{s_{0}^{\prime },\dots ,s_{n}^{\prime }\}=\{1,\dots ,2n+2\}\setminus
\{s_{0},\dots ,s_{n}\}$, $s_{0}^{\prime }<\dots <s_{n}^{\prime }$.
We have $%
A_{I}\in \mathcal{R}\Big(\big\{a_{j}\big\}_{j=1}^{2n+2}\Big)$ and $A_{I}%
\not\equiv%
%
0$ , since $\{a_{j}\}_{j=1}^{2n+2}$ are in general position.

Set $L=\{I\subset Q,\#I=n+1\}$, then $\#L=N:=\left( 
\begin{array}{c}
2n+2 \\ 
n+1
\end{array}
\right) $.

By the Laplace expansion Theorem, we have 
\begin{equation}
\sum_{I\in L}A_{I}h_{I}\equiv 0.  \tag{12}
\end{equation}

Let $I$, $J$, $K$ be distinct in $L$ . It is easy to see that:\vskip0.15cm%

$\left( (I\cup J)\diagdown (I\cap J)\right) \cap \left( (J\cup K)\diagdown
(J\cap K)\right) \cap \left( (K\cup I)\diagdown (K\cap I)\right) =\phi $ .%
\\
So, $C_{IJ}\cup C_{JK}\cup C_{KI}=\{1,...,n+2\}$ , where $%
C_{IJ}=\{1,...,n+2\}\diagdown \left( (I\cup J)\diagdown (I\cap J)\right) .$%
\\
Since $\text{dim}\big({}^{\ast }E^{i}\cap {}^{\ast }E^{j}\big)\leq m-2$ for
all $i\neq j$, $i\in \{1,\dots ,n+3\}$, $j\in \{1,\dots ,3n+1\},$ ${}^{\ast
}E^{i}\in \big\{{}^{\ast }E_{f}^{i},{}^{\ast }E_{g}^{i}\big\}$, ${}^{\ast
}E^{j}\in \big\{{}^{\ast }E_{f}^{j},{}^{\ast }E_{g}^{j}\big\}$ , and $f=g$
on $\bigcup\limits_{i=1}^{n+2}\big({}^{\ast }E_{f}^{i}\cap {}^{\ast
}E_{g}^{i}\big)$ (note that in the case $n=1$ we also have $f=g$ on $\bigcup\limits_{i=1}^{3}%
\big({}^{\ast }E_{f}^{i}\cap {}^{\ast }E_{g}^{i}\big)$), we have : 
\begin{equation*}
N_{\frac{h_{I}}{h_{J}}-1}^{[1]}(r)+N_{\frac{h_{J}}{h_{K}}-1}^{[1]}(r)+N_{%
\frac{h_{K}}{h_{I}}-1}^{[1]}(r)+\sum_{j\in
Q}{}^{>M}N_{(f,a_{j})}^{[1]}(r)+%
\sum_{k=1}^{3n+1}{}^{>M}N_{(g,a_{k})}^{[1]}(r)
\end{equation*}
\begin{equation}
\geq \sum_{i=1}^{n+2}{}^{\leq M}N_{(f,a_{i})}^{[1]}(r)  \tag{13}
\end{equation}
\\
Indeed, for $i\in \{1,...,n+2\}$, we may assume that $i\in C_{IJ}$. Let $%
z_{0}\in {}^{\ast }E_{f}^{i}$. If $z_{0}$ is not taken into account by $%
\underset{j\in Q}{\sum }{}^{>M}N_{(f,a_{j})}^{[1]}(r)$ or by$\underset{k=1}{%
\overset{3n+1}{\sum }}{}^{>M}N_{(g,a_{k})}^{[1]}(r)$ (this means that $%
v_{(f,a_{j})}(z_{0})\leq M$ and $v_{(g,a_{k})}(z_{0})\leq M$ for all $j\in Q$
, $k\in \{1,...,3n+1\})$ then $z_{0}\in {}^{\ast }E_{g}^{i}$ and by
omitting an analytic set of codimension $\geq 2$, we may assume that $%
(f,a_{j})(z_{0})\neq 0$ and $(g,a_{k})(z_{0})\neq 0$ for all $j\in
Q\diagdown \{i\}$ , $k\in \{1,...,3n+1\}\diagdown \{i\}.$ In particular $%
(f,a_{j})(z_{0})\neq 0$ , $(g,a_{j})(z_{0})\neq 0$ for all $j\in (I\cup
J)\diagdown (I\cap J).$\\
On the other hand, $f(z_{0})=g(z_{0}).$ Hence, $\frac{h_{I}}{h_{J}}(z_{0})=1$%
, this means that $z_{0}$ is taken into account by $N_{\frac{h_{I}}{h_{J}}%
-1}^{[1]}(r)$ , so we get (13).\\
By Lemma 3.2 and the First Main Theorem, we have : 
\begin{eqnarray*}
T_{f}(r) &\leq &\sum_{i=1}^{n+2}{}N_{(f,a_{i})}^{[n]}(r)+o\left(
T_{f}(r)\right)  \\
&\leq &\frac{M}{M+1}\sum_{i=1}^{n+2}{}^{\leq M}N_{(f,a_{i})}^{[n]}(r)+\frac{n%
}{M+1}\sum_{i=1}^{n+2}{}N_{(f,a_{i})}(r)+o\left( T_{f}(r)\right)  \\
&\leq &\frac{Mn}{M+1}\sum_{i=1}^{n+2}{}^{\leq M}N_{(f,a_{i})}^{[1]}(r)+\frac{%
n(n+2)}{M+1}T_{f}(r)+o\left( T_{f}(r)\right) .
\end{eqnarray*}
\\

Thus, we have
\begin{equation}
\frac{M+1-n(n+2)}{nM}T_{f}(r)\leq \sum_{i=1}^{n+2}{}^{\leq
M}N_{(f,a_{i})}^{[1]}(r)+\circ \big(T_{f}(r)\big).  \tag{14}
\end{equation}
\\
By (13) and (14) we have : 
\begin{equation*}
N_{\frac{h_{I}}{h_{J}}-1}^{[1]}(r)+N_{\frac{h_{J}}{h_{K}}-1}^{[1]}(r)+N_{%
\frac{h_{K}}{h_{I}}-1}^{[1]}(r)\geq \frac{M+1-n(n+2)}{nM}T_{f}(r)-\sum_{j\in
Q}{}^{>M}N_{(f,a_{j})}^{[1]}(r)
\end{equation*}
\begin{equation}
-\sum_{k=1}^{3n+1}{}^{>M}N_{(g,a_{k})}^{[1]}(r)-o\left( T_{f}(r)\right) . 
\tag{15}
\end{equation}
\\
We introduce an equivalence relation on $L$: $I\backsim J$ if and only if $%
\frac{h_{I}}{h_{J}}\in \mathcal{R}\Big(\big\{a_{j}\big\}_{j=1}^{3n+1}\Big).$%
\\
Set $\{L_{1},...,L_{s}\}=L\diagup _{\backsim }$ , ( $s\leq N:=\left( 
\begin{array}{c}
2n+2 \\ 
n+1
\end{array}
\right) $).

In order to prove (11), we show that $\#L_{v}\geq 2$ for all $v\in \{1,\dots
,s\}$. For each $v\in \{1,\dots ,s\}$, choose $I_{v}\in L_{v}$ and set 
\begin{equation*}
\sum_{I\in _{v}}A_{I}h_{I}=B_{v}h_{I_{v}},\quad B_{v}\in \mathcal{R}\Big(%
\big\{a_{j}\big\}_{j=1}^{3n+1}\Big).
\end{equation*}
Then (12) can be written as 
\begin{equation}
\sum_{v=1}^{s}B_{v}h_{I_{v}}\equiv 0.  \tag{16}
\end{equation}

+) If $B_{v}\equiv 0$ for all $v\in \{1,\dots ,s\}$, then $\#L_{v}\geq 2$
for all $v\in \{1,\dots ,s\}$ by $A_{I}%
\not\equiv%
%
0$, $I\in L$. We get (11).

+) If there exists some $B_{v}%
\not\equiv%
%
0$, then by (16) there are at least 3 of the $B_{1},\dots ,B_{s}$ different
from zero since $h_{I}%
\not\equiv%
%
0$ , $\dfrac{h_{I_{i}}}{h_{I_{j}}}\notin \mathcal{R}\Big(\big\{a_{j}\big\}%
_{j=1}^{3n+1}\Big)$, ($1\leq i\neq j\leq s,$ $I\in L$ ). \\
We want to apply Lemma 3.1 to (16), without loss of generality we may assume
that $B_{v}%
\not\equiv%
%
0$ for all $v\in \{1,\dots ,s\}$.

For each $\{i,j,k\}\subset \{1,\dots ,s\}$, set 
\begin{equation*}
T(r)=T_{\frac{B_{i}}{B_{j}}}(r)+T_{\frac{B_{j}}{B_{k}}}(r)+T_{\frac{B_{k}}{%
B_{i}}}(r)
\end{equation*}
then $T(r)=\circ \big(T_{f}(r)\big)$ as $r\rightarrow \infty $.

It is clear that $\dfrac{h_{I_{i}}}{h_{I_{j}}}-1%
\not\equiv%
%
0$, $\dfrac{h_{I_{j}}}{h_{I_{k}}}-1%
\not\equiv%
%
0$, $\dfrac{h_{I_{k}}}{h_{I_{i}}}-1%
\not\equiv%
%
0$.\\
By (10), (15), Theorem 5.2.29 in [7] and the First Main Theorem, we have: 
\begin{align}
& 3\cdot T_{[B_{i}h_{I_{i}}:B_{j}h_{I_{j}}:B_{k}h_{I_{k}}]}(r)+0(1)\geq T_{%
\frac{B_{i}h_{I_{i}}}{B_{j}h_{I_{j}}}}(r)+T_{\frac{B_{j}h_{I_{j}}}{%
B_{k}h_{I_{k}}}}(r)+T_{\frac{B_{k}h_{I_{k}}}{B_{i}h_{I_{i}}}}(r)  \notag \\
& \geq T_{\frac{h_{I_{i}}}{h_{I_{j}}}}(r)+T_{\frac{h_{I_{j}}}{h_{I_{k}}}%
}(r)+T_{\frac{h_{I_{k}}}{h_{I_{i}}}}(r)-T(r)  \notag \\
& \geq N_{\frac{h_{I_{i}}}{h_{i_{j}}}-1}(r)+N_{\frac{h_{I_{j}}}{h_{I_{k}}}%
-1}(r)+N_{\frac{h_{I_{k}}}{h_{I_{i}}}-1}(r)-\circ \big(T_{f}(r)\big)  \notag
\\
& \overset{(15)}{\geq }\frac{M+1-n(n+2)}{nM}T_{f}(r)-\sum_{j\in
Q}{}^{>M}N_{(f,a_{j})}^{[1]}(r)-%
\sum_{j=1}^{3n+1}{}^{>M}N_{(g,a_{j})}^{[1]}(r)-\circ \big(T_{f}(r)\big) 
\notag \\
& \geq \frac{M+1-n(n+2)}{nM}T_{f}(r)-\frac{1}{M+1}\sum_{j\in
Q}N_{(f,a_{j})}(r)-\frac{1}{M+1}\sum_{j=1}^{3n+1}N_{(g,a_{j})}(r)-\circ \big(%
T_{f}(r)\big)  \notag
\end{align}
\begin{align}
& \geq \frac{M+1-n(n+2)}{n\cdot M}T_{f}(r)-\frac{2(n+1)}{M+1}T_{f}(r)-\frac{%
3n+1}{M+1}T_{g}(r)-\circ \big(T_{f}(r)\big)  \notag \\
& \geq \Big(\frac{M+1-n(n+2)}{nM}-\frac{5n+3}{M+1}\Big)T_{f}(r)-\circ \big(%
T_{f}(r)\big),\ r\in A.  \tag{17}
\end{align}
Since $v_{(f,a_{j})}=v_{(g,a_{j})}$ on $E_{f}^{j}\cap E_{g}^{j}$, $j=1,\dots
,3n+1$, we have 
\begin{equation*}
\{z\in \mathbb{C}^{m}:h_{I}(z)=0\ \text{or}\ h_{I}(z)=\infty \}\subset
\bigcup_{j\in I}\{z\in \mathbb{C}^{m}:v_{(f,a_{j})}>M\ \text{or}\
v_{(g,a_{j})}(z)>M\}
\end{equation*}
for all $I\subset \{1,\dots ,3n+1\}$, $\#I=n+1$. Thus, we get
\begin{align*}
N_{h_{I}}^{[1]}(r)+N_{\frac{1}{h_{I}}}^{[1]}(r)& \leq \sum_{j\in
I}{}^{>M}N_{(f,a_{j})}^{[1]}(r)+\sum_{j\in I}{}^{>M}N_{(g,a_{j})}^{[1]}(r) \\
& \leq \frac{1}{M+1}\Big(\sum_{j\in I}N_{(f,a_{j})}(r)+\sum_{j\in
I}N_{(g,a_{j})}(r)\Big) \\
& \leq \frac{n+1}{M+1}\Big(T_{f}(r)+T_{g}(r)\Big)+0(1) \\
& \overset{(10)}{\leq }\frac{2(n+1)}{M+1}T_{f}(r)+0(1),\quad r\in A.
\end{align*}
So, we have
\begin{align}
& \sum_{v=1}^{s}N_{B_{v}h_{I_{v}}}^{[1]}(r)+\sum_{v=1}^{s}N_{\frac{1}{%
B_{v}h_{I_{v}}}}^{[1]}(r)\leq
\sum_{v=1}^{s}N_{h_{I_{v}}}^{[1]}(r)+\sum_{v=1}^{s}N_{\frac{1}{h_{I_{v}}}%
}^{[1]}(r)+  \notag \\
& +\sum_{v=1}^{s}N_{B_{v}}^{[1]}(r)+\sum_{v=1}^{s}N_{\frac{1}{B_{v}}%
}^{[1]}(r)\leq \sum_{v=1}^{s}\Big(N_{h_{I_{v}}}^{[1]}(r)+N_{\frac{1}{%
h_{I_{v}}}}^{1]}(r)\Big)+\circ \big(T_{f}(r)\big)  \notag \\
& \leq \frac{2s(n+1)}{M+1}T_{f}(r)+\circ \big(T_{f}(r)\big)\leq \frac{2(n+1)N%
}{M+1}T_{f}(r)+\circ \big(T_{f}(r)\big),\quad r\in A.  \tag{18}
\end{align}
By (17), (18) we have 
\begin{align}
& \sum_{v=1}^{s}N_{B_{v}h_{v}}^{[1]}(r)+\sum_{v=1}^{s}N_{\frac{1}{B_{v}h_{v}}%
}^{[1]}(r)\leq \frac{6n(n+1)NM}{(M+1)^{2}-n(n+2)(M+1)-n(5n+3)M}\cdot  \notag
\\
& \cdot T_{[B_{i}h_{I_{i}}:B_{h}h_{I_{j}}:B_{k}h_{I_{k}}]}(r)+\circ \big(%
T_{[B_{i}h_{I_{i}}:B_{j}h_{I_{j}}:B_{k}h_{I_{k}}]}(r)\big)  \notag \\
& <\frac{1}{N(N-1)}T_{\varphi _{ijk}}(r)\leq \frac{1}{s(s-1)}T_{\varphi
_{ijk}}(r),\quad r\in A,  \tag{19}
\end{align}
where \ $\varphi _{ijk}:=[B_{i}h_{I_{i}}:B_{j}h_{I_{j}}:B_{k}h_{I_{k}}]$.

Then by applying Lemma 3.1 to (16) we get: For each $i\in \{1,\dots ,s\}$
there exists $j\in \{1,\dots ,s\}$, $j\neq i$ such that $\frac{B_{i}h_{I_{i}}%
}{B_{j}h_{I_{j}}}$ is constant.

So $\frac{h_{I_{i}}}{h_{I_{j}}}\in \mathcal{R}\Big(\big\{a_{j}\big\}%
_{j=1}^{3n+1}\Big)$, this means that $L_{i}\cap L_{j}\neq \emptyset $. This
is a contradiction.

We have completed proof of (11). $\square $\\

Let $%
\mathcal{M}%
%
^{\ast }$ be the abelian multiplication group of all nonzero meromorphic
functions on $\mathbb{C}^{m}$. Define $%
\mathcal{H}%
%
$ $\subset 
\mathcal{M}%
%
^{\ast }$ by the set of all $h\in 
\mathcal{M}%
%
^{\ast }$ with $h^{k}\in \mathcal{R}\Big(\big\{a_{j}\big\}_{j=1}^{3n+1}\Big)$
for some positive integer $k$. It is easy to see that $%
\mathcal{H}%
%
$ is a subgroup of $%
\mathcal{M}%
%
^{\ast }.$\\
We have 
\begin{equation*}
\mathcal{M}%
%
^{\ast }\cap \mathcal{R}\Big(\big\{a_{j}\big\}_{j=1}^{3n+1}\Big)\subset 
\mathcal{H}%
%
\subset \widetilde{\mathcal{R}}\Big(\big\{a_{j}\big\}_{j=1}^{3n+1}\Big),
\end{equation*}
and the multiplication group $%
\mathcal{G}%
%
:=%
\mathcal{M}%
%
^{\ast }\diagup 
\mathcal{H}%
%
$ is a torsion free abebian group. We denote by $[h]$ the class in $%
\mathcal{G}%
%
$ containing $h\in 
\mathcal{M}%
%
^{\ast }$. Consider the subgroup $\widetilde{%
\mathcal{G}%
%
}$ of $%
\mathcal{G}%
%
$ generated by $[h_{1}],\dots ,[h_{3n+1}]$ and choose suitable functions $%
\eta _{1},\dots ,\eta _{t}\in 
\mathcal{M}%
%
^{\ast }$ such that $[\eta _{1}],\dots ,[\eta _{t}]$ give a basis of $%
\widetilde{%
\mathcal{G}%
%
}$. Then each $h_{j}$ can be uniquely represented as $h_{j}=c_{j}\eta
_{1}^{\ell _{j_{1}}}\cdots \eta _{t}^{\ell _{j_{t}}}$, $c_{j}\in 
\mathcal{H}%
%
$, $\ell _{j_{r}}\in \mathbb{Z}$. For these integers $\ell _{j_{r}}$ we can
choose suitable integers $p_{1},\dots ,p_{t}$ satisfying the condition: For
integers $\ell _{j}=p_{1}\ell _{j_{1}}+\dots +p_{t}\ell _{j_{t}}$, $(1\leq
j\leq 3n+1)$, $\ell _{i}=\ell _{j}$ if and only if $(\ell _{i_{1}},\dots
,\ell _{i_{t}})=(\ell _{j_{1}},\dots ,\ell _{j_{t}})$, or equivalently 
\begin{equation*}
\frac{h_{i}}{h_{j}}\in 
\mathcal{H}%
%
.
\end{equation*}

We now show that:\\
There is a subset $I_{0}=\{j_{0},...,j_{n}\}\subset $ $Q_{0}$ such that
\begin{equation}
\frac{h_{i}}{h_{j}}\in 
\mathcal{H}%
%
\quad \text{for all}\quad i,j\in I_{0}.  \tag{20}
\end{equation}
We assume that, after a suitable change of indices, we have $\ell _{1}\leq \dots \leq
\ell _{3n+1}$.

Take the subset $Q=\{1,\dots ,n+1,2n+1,\dots ,3n+1\}$ of $Q_{0}$ which
contains $(2n+2)$ elements and apply (11) to the $h_{j}^{\prime }$ $(j\in Q)
$ to show that there is a subset $\{i_{0},\dots ,i_{n}\}$ of $Q$ satisfying
the condition that $\{i_{0},\dots ,i_{n}\}\neq \{1,\dots ,n+1\}$, $%
i_{0}<\dots <i_{n}$ and 
\begin{equation*}
\frac{h_{i_{0}}\cdots h_{i_{n}}}{h_{1}\cdots h_{n+1}}\in \mathcal{R}\Big(%
\big\{a_{j}\big\}_{j=1}^{3n+1}\Big).
\end{equation*}
From this, it follows that
\begin{equation*}
(\ell _{i_{0}}-\ell _{1})+\dots +(\ell _{i_{n}}-\ell
_{n+1})=\sum_{s=1}^{t}p_{s}(\ell _{i_{0_{s}}}+\dots +\ell _{i_{n_{s}}}-\ell
_{1_{s}}-\dots -\ell _{n+1_{s}})=0.
\end{equation*}
Since $\ell _{i_{0}}\geq \ell _{1},\dots ,\ell _{i_{n}}\geq \ell _{n+1}$ ,
this is posible only if $\ell _{n+1}=\ell _{i_{n}}$ so $\ell _{n+1}=\dots
=\ell _{2n+1}$ (note that $i_{n}\geq 2n+1$ ). Then take $I_{0}=%
\{n+1,...,2n+1\}$, so we get (20). $\square $

Let $u_{i}=\dfrac{h_{j_{i}}}{h_{j_{0}}}\in 
\mathcal{H}%
%
\subset \widetilde{\mathcal{R}}\Big(\big\{a_{j}\big\}_{j=1}^{3n+1}\Big)$, $%
i\in \{0,\dots ,n\}$. Then we have:

\begin{equation}
\left\{ 
\begin{array}{c}
a_{j_{i}0}f_{0}+\dots +a_{j_{i}n}f_{n}=u_{i}h_{j_{0}}(a_{j_{i}0}g_{0}+\dots
+a_{j_{i}n}g_{n}) \\ 
i=0,\dots ,n.
\end{array}
\right.  \tag{21}
\end{equation}

+)Assume $n=1.$

Set $A=\left( 
\begin{array}{ll}
a_{j_{0}0} & a_{j_{0}1} \\ 
a_{j_{1}0} & a_{j_{1}1}
\end{array}
\right) ,B=\left( 
\begin{array}{ll}
u_{0} & 0 \\ 
0 & u_{1}
\end{array}
\right) $.\\

By (21) , $A.\left( 
\begin{array}{l}
f_{0} \\ 
f_{1}
\end{array}
\right) =h_{j_{0}}.B.A.\left( 
\begin{array}{l}
g_{0} \\ 
g_{1}
\end{array}
\right) \Rightarrow A^{-1}.B^{-1}.A.\left( 
\begin{array}{l}
f_{0} \\ 
f_{1}
\end{array}
\right) =h_{j_{0}}.\left( 
\begin{array}{l}
g_{0} \\ 
g_{1}
\end{array}
\right) $\\
We get 1) of the Theorem 1 (with $L=A^{-1}.B^{-1}.A).$

+) Assume $n\geq 2.$\\
Set $F=\left( (f,\widetilde{a_{j_{0}}}):...:(f,\widetilde{a_{j_{n}}})\right) $ and $G=\left(
(g,\widetilde{a_{j_{0}}}):...:(g,\widetilde{a_{j_{n}}})\right) $. They are meromorphic mappings of $%
\mathbb{C}^{m}$ into $\mathbb{C}P^{n}$ . Take meromorphic functions $h,u$ on 
$\mathbb{C}^{m}$ such that $F=\left( \frac{(f,\widetilde{a_{j_{0}}})}{h}:...:\frac{%
(f,\widetilde{a_{j_{n}}})}{h}\right) ,G=\left( \frac{(g,\widetilde{a_{j_{0}}})}{u}:...:\frac{%
(g,\widetilde{a_{j_{n}}})}{u}\right) $ are reduced representations. By Lemma 3.2 we have $N_{h}(r)=o\left( T_{f}(r)\right) $,
 $N_{u}(r)=o\left( T_{f}(r)\right) $, $N_{\frac{1}{h}}(r)=o(T_f(r))$ and $N_{\frac{1}{u}}(r)=o(T_f(r))$.\\
Put $F_{i}:=\frac{(f,\widetilde{a_{j_{i}}})}{h},G_{i}:=\frac{(g,\widetilde{a_{j_{i}}})}{u},i\in
\{0,...,n\}.$\\
Since $u_{i}\in 
\mathcal{H}%
%
$, $i\in \{0,\dots ,n\}$, we can choose a positive $k$ such that $%
(u_{i})^{k}\in \mathcal{R}\Big(\big\{a_{j}\big\}_{j=1}^{3n+1}\Big)$ for all $%
i\in \{0,...,n)$ .\\
By (21) \ we have 
\begin{equation*}
\left\{ 
\begin{array}{c}
F_{i}=\frac{u_{i}h_{j_{0}}uG_{i}}{h} \\ 
i=0,...,n
\end{array}
\right.
\end{equation*}
\\
Since $G=(G_{0}:...:G_{n})$ is a reduced representation and $F_{i}$ $%
(i=1,...,n)$ are holomorphic functions , we have: 
\begin{equation*}
N_{\frac{1}{h_{j_{0}}}}(r)\leq \sum_{i=0}^{n}N_{u_{i}}(r)+N_{u}(r)+N_{\frac{1}{h}}(r)\leq
\sum_{i=0}^{n}N_{(u_{i})^{k}}(r)+N_{u}(r)+N_{\frac{1}{h}}(r)
\end{equation*}
\begin{equation}
\leq \sum_{i=0}^{n}T_{(u_{i})^{k}}(r)+0(1)+N_{u}(r)+N_{\frac{1}{h}}(r)=o\left( T_{f}(r)\right). 
\tag{22}
\end{equation}
\\
Suppose that $F%
\not\equiv%
%
G$ , then there exist $0\leq s<v\leq n$ such that :

$\left| 
\begin{array}{ll}
F_{s} & F_{v} \\ 
G_{s} & G_{v}
\end{array}
\right| 
\not\equiv%
%
0\Rightarrow $ $\left( \frac{h}{uh_{j_{0}}u_{v}}-\frac{h}{uh_{j_{0}}u_{s}}%
\right) F_{s}F_{v}%
\not\equiv%
%
0\,$ .\\
Define the meromorphic mapping $F\wedge G:=$ $(....:\left| 
\begin{array}{ll}
F_{i} & F_{j} \\ 
G_{i} & G_{j}
\end{array}
\right| :...):$ $\mathbb{C}^{m}\longrightarrow $ $\mathbb{C}P^{N_{2}}$, ( $%
0\leq i<j\leq n$ , $N_{2}=\left( 
\begin{array}{c}
n+1 \\ 
2
\end{array}
\right) -1$ ) .\\
Take $\mu _{F\wedge G}$ a holomophic function on $\mathbb{C}^{m}$ such that $%
(....:\frac{1}{\mu _{F\wedge G}}\left| 
\begin{array}{ll}
F_{i} & F_{j} \\ 
G_{i} & G_{j}
\end{array}
\right| :...)$ is a reduded representation of $F\wedge G$ .\\
It is easy to see that there exists a subset $I_{sv}\subseteq
\{1,...,n+4\}\diagdown \{j_{s},j_{v}\}$ such that 
\begin{equation*}
\#I_{sv}=n+2,\#\left( \{1,...,n+3\}\diagdown (I_{sv}\cup \{j_{s}\})\ \
\right) \leq 1,\text{ and}
\end{equation*}

\begin{equation}
\#\left( \{1,...,n+3\}\diagdown (I_{sv}\cup \{j_{v}\})\ \
\right) \leq 1.  \tag{23}
\end{equation}

In fact, we take $I_{sv}=\{1,...,n+2\}$ if $\{j_{s},j_{v}\}\cap
\{1,...,n+3\}=\phi ,$\\
$I_{sv}=\{1,...,n+3\}\diagdown \{j_{s},j_{v}\}$ if$\ \#\left(
\{j_{s},j_{v}\}\cap \{1,...,n+3\}\right) =1$ and

$I_{sv}=\{1,...,n+4\}\diagdown \{j_{s},j_{v}\}$ if $\{j_{s},j_{v}\}\subset
\{1,...,n+3\}.$\\
By assumptions ii) and iii) we have :

\begin{equation}
N_{\frac{1}{F_{s}F_{v}}\mu _{F\wedge G}}(r)\geq \sum_{i\in I_{sv}}{}^{\leq
M}N_{(f,a_{i})}^{[1]}(r)-\sum_{j\in
\{j_{s},j_{v}\}}{}^{>M}N_{(f,a_{j})}^{[1]}(r)-\sum_{i\in
I_{sv}}{}^{>M}N_{(g,a_{i})}^{[1]}(r)  \tag{24}
\end{equation}
\\

Indeed, for $i_{0}\in I_{sv}$, let $z_{0}\in {}^{\ast }E_{f}^{i_{0}}$ be a
generic point of a component $D$ of ${}^{\ast }E_{f}^{i_{0}}$. If $z_{0}$ is
not taken into account by $\underset{j\in \{j_{s},j_{v}\}}{\sum }%
{}^{>M}N_{(f,a_{j})}^{[1]}(r)$ or by $\underset{i\in I_{sv}}{\overset{}{\sum 
}}{}^{>M}N_{(g,a_{i})}^{[1]}(r)$ (this means that $v_{(f,a_{j})}(z_{0})\leq M
$, $j\in \{j_{s},j_{v}\}$ and $\ v_{(g,a_{i})}(z_{0})\leq M$, $i\in I_{sv})$
then $z_{0}\in {}^{\ast }E_{g}^{i_{0}}$ (which implies $f(z_{0})=g(z_{0})$).
Since $z_{0}\in D$ is generic, we can omit an analytic set of codimension $%
\geq 2$, so we may assume that $(f,a_{j_{s}})(z_{0})\neq 0$ ,$%
(f,a_{j_{v}})(z_{0})\neq 0$ (note that by (23) we cannot have $%
\{i_{0},j_{s}\}\subset \{n+4,...,3n+1\}$ or $\{i_{0},j_{v}\}\subset
\{n+4,...,3n+1\}$), which implies $F_{s}(z_{0})\neq 0$ , $F_{v}(z_{0})\neq 0$
. Since we have $f(z_{0})=g(z_{0})$ on $D$, we get $\mu _{F\wedge G}($ $%
z_{0})=0$ on $D$. This means that $z_{0}$ is taken into account by $N_{\frac{%
1}{F_{s}F_{v}}\mu _{F\wedge G}}(r)$ , so we get (24).\\
So, we have
\begin{eqnarray*}
N_{\frac{1}{F_{s}F_{v}}\mu _{F\wedge G}}(r) &\geq &\sum_{i\in
I_{sv}}{}^{\leq M}N_{(f,a_{i})}^{[1]}(r)-\frac{1}{M+1}\left( \sum_{j\in
\{j_{s},j_{v}\}}{}N_{(f,a_{j})}(r)-\sum_{i\in
I_{sv}}{}N_{(g,a_{i})}(r)\right)  \\
&\geq &\sum_{i\in I_{sv}}{}^{\leq M}N_{(f,a_{i})}^{[1]}(r)-\frac{2}{M+1}{}%
T_{f}(r)-\frac{n+2}{M+1}T_{g}(r)
\end{eqnarray*}
\begin{equation}
\geq \sum_{i\in I_{sv}}{}^{\leq M}N_{(f,a_{j})}^{[1]}(r)-\frac{(n+4)}{M}{}%
T_{f}(r),r\in A.  \tag{25}
\end{equation}
\\
By Lemma 3.2 and the First Main Theorem, we have:\\
\begin{eqnarray*}
T_{f}(r) &\leq &\sum_{i\in I_{sv}}{}N_{(f,a_{i})}^{[n]}(r)+o\left(
T_{f}(r)\right)  \\
&\leq &\frac{M}{M+1}\sum_{i\in I_{sv}}{}^{\leq M}N_{(f,a_{i})}^{[n]}(r)+%
\frac{n}{M+1}\sum_{i\in I_{sv}}{}N_{(f,a_{i})}(r)+o\left( T_{f}(r)\right)  \\
&\leq &\frac{Mn}{M+1}\sum_{i\in I_{sv}}{}^{\leq M}N_{(f,a_{i})}^{[1]}(r)+%
\frac{n(n+2)}{M+1}T_{f}(r)+o\left( T_{f}(r)\right) .
\end{eqnarray*}
\\
\begin{equation*}
\Rightarrow \left( \frac{M+1-n(n+2)}{Mn}\right) T_{f}(r)\leq \sum_{i\in
I_{sv}}{}^{\leq M}N_{(f,a_{i})}^{[1]}(r)+o\left( T_{f}(r)\right) .
\end{equation*}
So, by (25) we have: 
\begin{equation}
\left( \frac{M+1-2n(n+3)}{Mn}\right) T_{f}(r)\leq N_{\frac{1}{F_{s}F_{v}}\mu
_{F\wedge G}}(r)+o\left( T_{f}(r)\right) ,r\in A.  \tag{26}
\end{equation}
\\
By the definition of $\mu _{F\wedge G}$ , we have:

\begin{eqnarray*}
N_{\frac{1}{F_{s}F_{v}}\mu _{F\wedge G}}(r) &\leq &N_{\frac{1}{F_{s}F_{v}}%
\left| 
\begin{array}{ll}
F_{s} & F_{v} \\ 
G_{s} & G_{v}
\end{array}
\right| }(r)=N_{\left( \frac{h}{uh_{j_{0}}u_{v}}-\frac{h}{uh_{j_{0}}u_{s}}%
\right) }(r) +N_{\frac{1}{u}}(r)\\
&\leq &N_{\left( \frac{1}{u_{v}}-\frac{1}{u_{s}}\right) }(r)+N_{\frac{1}{%
h_{j_{0}}}}(r)+N_{h}(r) +N_{\frac{1}{u}}(r)\\
&\overset{(22)}{\leq } &N_{\left( \frac{1}{(u_{v})^{k}}-\frac{1}{(u_{s})^{k}}%
\right) }(r)+o\left( T_{f}(r)\right) \leq T_{\left( \frac{1}{(u_{v})^{k}}-%
\frac{1}{(u_{s})^{k}}\right) }(r)+o\left( T_{f}(r)\right) \\
&\leq &T_{\frac{1}{(u_{v})^{k}}}(r)+T_{\frac{1}{(u_{s})^{k}}}(r)+o\left(
T_{f}(r)\right) =o\left( T_{f}(r)\right) .
\end{eqnarray*}
\\

This contracdicts to (26). Thus $F=G$ $\Rightarrow $ $f=g$ , so we get 2) of
Theorem 1$.$ This completes the proof of Theorem 1. $\square $\\

We can obtain Theorem 2 by an argument similar to the proof of Theorem 1
with the following remarks:

+) We do not need the assumption (10).

+) Similarly to (13) we have : 
\begin{align*}
N_{\frac{h_{I}}{h_{J}}-1}^{[1]}(r)+N_{\frac{h_{J}}{h_{K}}-1}^{[1]}(r)& +N_{%
\frac{h_{K}}{h_{I}}}^{[1]}(r)+\sum_{j\in Q}{}^{>M}N_{(f,a_{j})}^{[1]}(r) \\
& \geq \sum_{i=1}^{n+2}{}^{\leq M}N_{(f,a_{i})}^{[1]}(r)\quad \text{for all }%
r.
\end{align*}
So, similarly to (17) we have, 
\begin{equation}
3\cdot T_{[B_{i}h_{I_{i}}:B_{j}h_{I_{j}}:B_{k}h_{I_{k}}]}(r)\geq \Big(\frac{%
(M+1)-n(n+2)}{nM}-\frac{2(n+1)}{M+1}\Big)T_{f}(r)+\circ \big(T_{f}(r)\big), 
\tag{17'}
\end{equation}

+) Since $\min \big\{v_{(f,a_{j})},M\big\}=\min \big\{v_{(g,a_{j})},M\big\}$%
, $j\in \{1,\dots ,3n+1\}$, we have 
\begin{equation*}
\big\{z\in \mathbb{C}^{m}:h_{I}(z)=0\ \text{or}\ h_{I}(z)=\infty \big\}%
\subset \bigcup_{j\in I}\big\{z\in \mathbb{C}^{m}:v_{(f,a_{j})}(z)>M\big\}
\end{equation*}
for all $I\subset \{1,\dots ,3n+1\}$, $\#I=n+1$.

Thus, we get
\begin{equation*}
N_{h_{I}}^{[1]}(r)+N_{\frac{1}{h_{I}}}^{[1]}(r)\leq \sum_{j\in
I}{}^{>M}N_{(f,a_{j})}^{[1]}(r)\leq \frac{1}{M+1}\sum_{j\in
I}N_{(f,a_{j})}(r)\leq \frac{n+1}{M+1}T_{f}(r)+0(1)
\end{equation*}
So, similarly to (18) we have, 
\begin{equation}
\sum_{v=1}^{s}N_{B_{v}h_{I_{v}}}^{[1]}(r)+\sum_{v=1}^{s}N_{\frac{1}{%
B_{v}h_{I_{v}}}}^{[1]}(r)\leq \frac{N(n+1)}{M+1}T_{f}(r)+\circ \big(T_{f}(r)%
\big)\quad \text{for all}\ r.  \tag{18'}
\end{equation}
By (17') and (18') we have : 
\begin{equation*}
\sum_{v=1}^{s}N_{B_{v}h_{I_{v}}}^{[1]}(r)+\sum_{v=1}^{s}N_{\frac{1}{%
B_{v}h_{I_{v}}}}^{[1]}(r)\leq \frac{1}{s(s-1)}T_{\varphi _{ijk}}(r)\quad 
\text{for all}\ r.
\end{equation*}

+) Similarly to (24) we have:\\
\begin{equation*}
N_{\frac{1}{F_{s}F_{v}}\mu _{F\wedge G}}(r)\geq \sum_{i\in I_{sv}}{}^{\leq
M}N_{(f,a_{i})}^{[1]}(r)-\sum_{j\in
\{j_{s},j_{v}\}}{}^{>M}N_{(f,a_{j})}^{[1]}(r).
\end{equation*}
So, similarly to (25) we have: 
\begin{equation*}
N_{\frac{1}{F_{s}F_{v}}\mu _{F\wedge G}}(r)\geq \sum_{i\in I_{sv}}{}^{\leq
M}N_{(f,a_{i})}^{[1]}(r)-\frac{2}{M}{}T_{f}(r).
\end{equation*}

 \noindent Gerd Dethloff, Tran Van Tan \\
Universit\'{e} de Bretagne Occidentale \\
  UFR Sciences et Techniques \\
D\'{e}partement de Math\'{e}matiques \\
6, avenue Le Gorgeu, BP 452 \\
  29275 Brest Cedex, France \\
e-mail: gerd.dethloff@univ-brest.fr\\

\end{document}